\documentclass[12pt]{article}
    \usepackage{t1enc}
    \usepackage[latin1]{inputenc}
    \usepackage[english]{babel}
        \usepackage{latexsym}
        \newcommand{\bp}{{\bf P}}

\newcommand{\be}{{\bf E}}

\newtheorem{theorem}{Theorem}[section]

\newtheorem{lemma}{Lemma}[section]

\def\begg{begin{equation}}
\def\endd{\end{equation}}

\def\z2{{\cal Z}_2}

\def\begg{\begin{equation}}
\def\endd{\end{equation}}
\def\bege{\begin{eqnarray}}
\def\ende{\end{eqnarray}}

\def\ep{{\varepsilon}}
\def\al{{\alpha}}

\begin{document}

 \centerline{\Large\bf TRANSIENT NN RANDOM WALK ON
THE LINE}
\medskip

\bigskip \bigskip \bigskip \bigskip \bigskip

\renewcommand{\thefootnote}{1} \noindent
\textbf{Endre Cs\'{a}ki}\footnote{Research supported by the
Hungarian National Foundation for Scientif\/ic Research, Grant No.
  K 61052 and K 67961.}\newline Alfr\'ed R\'enyi Institute of
Mathematics, Hungarian Academy of Sciences, Budapest, P.O.B. 127,
H-1364, Hungary. E-mail address: csaki@renyi.hu

\bigskip

\renewcommand{\thefootnote}{2} \noindent \textbf{Ant\'{o}nia
F\"{o}ldes}\footnote{Research supported by a PSC CUNY Grant, No.
68030-0037.}\newline Department of Mathematics, College of Staten
Island, CUNY, 2800 Victory Blvd., Staten Island, New York 10314,
U.S.A. E-mail address: foldes@mail.csi.cuny.edu

\bigskip

\noindent \textbf{P\'al R\'ev\'esz}$^1$ \newline Institut f\"ur
Statistik und Wahrscheinlichkeitstheorie, Technische Universit\"at
Wien, Wiedner Hauptstrasse 8-10/107 A-1040 Vienna, Austria. E-mail
address: reveszp@renyi.hu

\bigskip \bigskip \bigskip

\noindent \textit{Abstract:} We prove strong theorems for the
local time at infinity of a nearest neighbor  transient random
walk. First, laws of the iterated logarithm are given for the
large values of the local time. Then we  investigate the length of
intervals over which the walk runs through (always from left  to
right) without ever  returning.
\bigskip


\bigskip

 \vspace{.1cm}

\noindent Running head: Transient NN random walk.  \vfill
\renewcommand{\thesection}{\arabic{section}.}
\section{Introduction}

\renewcommand{\thesection}{\arabic{section}} \setcounter{equation}{0}
\setcounter{theorem}{0} \setcounter{lemma}{0}

Let $X_0=0,\ X_1,X_2,\ldots$ be a Markov chain with
\begin{eqnarray}\label{def}
E_i&:=&\bp(X_{n+1}=i+1\mid X_n=i)=1-\bp(X_{n+1}=i-1\mid X_n=i)\\
&=&\left\{\begin{array}{ll} 1\quad & if\quad  i=0\\
\nonumber
 1/2+p_i\quad & if\quad i=1,2,\ldots,
\end{array}\right.
\end{eqnarray}

\noindent where $0\leq p_i< 1/2,\ i=1,2,\ldots$.  This sequence
$\{X_i\}$ describes the motion of a particle which starts at zero,
moves over the nonnegative integers and going away from 0 with a
larger probability than to the direction of 0. We will be
interested in the case when $\{p_i, \,i=1,2...\}$ goes to zero.
That is to say 0 has a repelling power which becomes small if the
particle is far away from 0. We intend to characterize the local
time of this motion.

A slightly different but symmetric variation of the same motion
can be defined  as follows.
 Let $X^*_0=0,\ X^*_1,X^*_2,\ldots$ be a Markov chain
with
\begin{eqnarray*}
E^*_i&:=&\bp(X^*_{n+1}=i+1\mid X^*_n=i)=1-\bp(X^*_{n+1}=i-1\mid
X^*_n=i)=\\
&=&\left\{\begin{array}{ll}
1/2\quad & if\quad  i=0,\\
1/2+p_i\quad & if\quad i=1,2,\ldots,\\
1/2-p_i\quad & if\quad i=-1,-2,\ldots
\end{array}\right.
\end{eqnarray*}
Our results can be rephrased  with the obvious modification for
this walk as well. However to be in line with the existing
literature we will use the definition in  (\ref{def}).

The properties of this model, often called birth and death chain,
connections with orthogonal polynomials in particular, has been treated
extensively in the literature. See e.g. the classical paper by Karlin and
McGregor \cite{KMG}, or more recent papers by Coolen-Schrijner and Van
Doorn \cite{C-SD} and Dette \cite{DE01}.

As it will turn out in this paper, the properties of the
walk and its local time is very sensitive even for small changes
in $\{p_i\}$-s. \noindent
There is a well-known  result in the literature (cf. e.g. Chung
\cite{CH}) characterizing those sequences $\{p_i\}$ for which $\{X_i\}$
is transient (resp. recurrent).

\smallskip

\noindent
 {\bf Theorem A:} (\cite{CH}, page 74) {\it Let $X_n$ be a Markov chain
 with transition probabilities given in {\rm (\ref{def}).}
Define

\begg U_i:={\frac{1-E_i}{E_i}}={\frac{1/2-p_i}{1/2+p_i}}
\label{uif}
\endd
Then $X_n$  is transient if and only if}
$$\sum_{k=1}^{\infty} \prod_{i=1}^k U_i < \infty.$$
This criteria however does not reveal explicitly what are  the
transient/reccurent type of $\{p_i\}$ sequences.
 \noindent Lamperti
\cite{LA60}, \cite{LA602} proved a more general theorem about
recurrence and transience of real nonnegative processes (not
necessarily Markov chains). Here we spell out his result in our
setup only, which easily follows from Theorem A as well.

\smallskip

\noindent  {\bf Corollary:} {\it If for all
 $i$ large enough,
\begg p_i\leq\frac{1}{4i}+O\left(\frac{1}{i^{1+\delta}}\right)
\qquad \delta>0,
\endd then $\{X_i\}$ is recurrent. If instead,  for some $\theta>1$
\begg p_i\geq \frac{\theta }{4i}
\endd
for $i$  large enough, then $\{X_i\}$ is transient.}

As we proceed to find the necessary tools for getting  results
 about the local time, as a byproduct,
 we will get a much sharper version of this Corollary.

In this paper we concentrate only on the transient case.

There are many results in the literature about the limiting behavior of
$\{X_n\}$, depending on the sequence $\{p_i\}.$ Lamperti \cite{LA63}
determined the limiting distribution of $X_n$.

\smallskip
\noindent {\bf Theorem B:} (\cite{LA63})
{\it If $\,\,\lim_{i\to\infty}\,ip_i=B/4>0$, then}
$$
\lim_{n\to \infty} \bp\left(\frac{X_n}{\sqrt{n}}<x\right)=
\frac{1}{2^{B/2-1/2}\Gamma(B/2+1/2)}\int_0^x u^Be^{-u^2/2}\, du.
$$

In fact, Lamperti \cite{LA63} (see also Rosenkrantz \cite{RO73}) proved
weak convergence of $X_n/\sqrt{n}$ to a Bessel process as well. We
intend to give further connections (strong invariance, etc.) between $X_n$
and Bessel process in a subsequent paper.

The law of the iterated logarithm for $X_n$ was given by Br\'ezis
et al. \cite{BRS}, Sz\'ekely \cite{SZ}, Gallardo \cite{GA84}, Voit
\cite{VO90}. Their somewhat more general results specialized in
our setup, reads as follows.


\smallskip

\noindent {\bf Theorem C:} (\cite{BRS}, \cite{SZ}, \cite{GA84},
\cite{VO90})

{\it If $\,\,\lim_{i\to\infty}\,ip_i=c>0$, then}
$$\limsup_{n\to \infty} \frac{X_n}{\sqrt{2n\log \log
n}}=1\qquad {\rm a.s.}$$

Voit \cite{VO92}  has proved a law of large numbers for certain
Markov chains, which we quote in our setup only.

\medskip
\noindent {\bf Theorem D:} (\cite{VO92})
{\it If $\,\,\lim_{i\to\infty}\,i^\alpha p_i=c>0$ for some $0<\alpha<1$,
then}
$$\lim_{n\to \infty}\frac{X_n}{n^{1/(1+\alpha)}}=
2c(1+\alpha)\qquad {\rm a.s.}$$

Our main concern in this paper is to study the local time of
$\{X_n\}$, defined by

\begg\xi(x,n):=\#\{k:0\leq k\leq n,\ X_k=x\},\quad x=0,1,2,\ldots,
\label{locn}\endd and \begg\xi(x,\infty):=\lim_{n\to
\infty}\xi(x,n). \label{loc1}
\endd

\noindent
\renewcommand{\thesection}{\arabic{section}.}
\section{Lemmas and Notations}

\renewcommand{\thesection}{\arabic{section}} \setcounter{equation}{0}
\setcounter{theorem}{0} \setcounter{lemma}{0}

For $U_i$ defined in (\ref{uif}) we get by elementary calculation
that

\smallskip

\noindent {\bf Fact 1.}
\begin{eqnarray}
 U_i&=&{\frac{1-E_i}{E_i}}={\frac{1/2-p_i}{1/2+p_i}}=1-4p_i+O(p_i^2)\nonumber
\\&=&\exp(-4p_i+O(p_i^2))\quad (i=0,\pm 1,\pm 2,\ldots).
\label{ui}
\end{eqnarray}

 Introduce the notation

$$D(m,n):=\left\{\begin{array}{ll}
0 & {\rm if\ } n=m,\\
1 & {\rm if\ } n=m+1,\\
1+\displaystyle{\sum_{j=1}^{n-m-1}\prod_{i=1}^j U_{m+i}}= & \\
1+\displaystyle{\sum_{j=1}^{n-m-1}}
\exp\left(-(1+o_m(1))4\sum_{i=m+1}^{m+j}p_i\right) & {\rm if\ }
n\geq m+2.
\end{array}\right.$$
Denote $$\lim_{n\to \infty}D(m,n)=:D(m,\infty).$$
\begin{lemma} If $p_i\downarrow 0$ then for $m$ large  enough
$$D(m,\infty) \geq \frac{C}{p_m}, $$
where $C$ is an absolute constant. Consequently,
 $$\lim_{m\to \infty}D(m,\infty)=
+\infty.$$
\end{lemma}

\noindent
{\bf Proof:} Let $\displaystyle{\frac{1}{p_m}\leq j\leq
\frac{2}{p_m}, }$ then from (\ref{ui}) for $m$ big enough we have
$$ \sum_{i=m}^{m+j}(p_i+C p_i^2)
\leq \frac{2}{p_m}(p_m+C p_m^2) \leq 2 ( 1+C p_m)\leq 2(1+C).
$$

\noindent Consequently

$$\exp\left(- \sum_{i=m}^{m+j}(p_i+C p_i^2)\right
)\geq \exp (-2(1+C))
$$
and
\begin{eqnarray*}
D(m,\infty)&\geq& 1+\sum_{j=0}^\infty \exp\left(
-\sum_{i=m}^{m+j}(p_i+C p_i^2)\right) \\&\geq&
1+\sum_{j=[\frac{1}{p_m}]} ^{[\frac{2}{p_m}]}\exp(-2(1+C))\geq
1+\frac{1}{p_m}\exp(-2(1+C)).
\end{eqnarray*}
$\Box$

\noindent

 For $0 \leq a\leq b\leq c$ define
\begin{eqnarray*}
\lefteqn{p(a,b,c):=}\\
&=&\bp(\min\{j:\ j>m,\ X_j=a\}<\min\{j:\ j>m,\ X_j=c\}\mid X_m=b),\\
&&
\end{eqnarray*}

\noindent i.e. $p(a,b,c)$ is the probability that a particle
starting from $b$ hits $a$ before $c$.

\begin{lemma} For $0 \leq a \leq b \leq c$

$$p(a,b,c)=1-{\frac{D(a,b)}{D(a,c)}}.$$

\noindent Especially

\begin{equation}
p(0,1,n)=1-{\frac{1}{D(0,n)}},\quad p(n,n+1,\infty)=
1-{\frac{1}{D(n,\infty)}}.\label{it0}
\end{equation}
\end{lemma}

\vspace{2ex}\noindent
{\bf Proof:} The proof of this lemma is fairly standard, we
 give it for completeness. Clearly we have

\begin{eqnarray*}
p(a,a,c)&=&1,\\
p(a,c,c)&=&0,\\
p(a,b,c)&=&E_bp(a,b+1,c)+(1-E_b)p(a,b-1,c).
\end{eqnarray*}

\noindent Consequently

$$p(a,b+1,c)={\frac{1}{E_b}}p(a,b,c)-{\frac{1-E_b}{E_b}}
p(a,b-1,c)$$

\noindent and
\begin{eqnarray*}
p(a,b+1,c)-p(a,b,c)&=&{\frac{1-E_b}{E_b}}(p(a,b,c)-p(a,b-1,c))=\\
&=&U_b(p(a,b,c)-p(a,b-1,c)).
\end{eqnarray*}

\noindent By iterat\/ion we get

\begin{eqnarray}
\lefteqn{p(a,b+1,c)-p(a,b,c)=} \label{it1}\\
&=&U_bU_{b-1}(p(a,b-1,c)-p(a,b-2,c))\nonumber\\
&=&\ldots=U_bU_{b-1}\cdots U_{a+1}(p(a,a+1,c)-p(a,a,c))=\nonumber\\
&=&U_bU_{b-1}\cdots U_{a+1}(p(a,a+1,c)-1).\nonumber
\end{eqnarray}

\noindent Starting with the trivial identity

$$
p(a,a+1,c)-p(a,a,c)=p(a,a+1,c)-1
$$

 \noindent and adding to it  the above equations for $b=a+1,\ldots,c-1$ we get

$$-1=p(a,c,c)-p(a,a,c)=D(a,c)(p(a,a+1,c)-1),$$

\noindent i.e.

\begin{equation}
p(a,a+1,c)=1-{\frac{1}{D(a,c)}}. \label{it2}
\end{equation}

\noindent Hence (\ref{it1}) and (\ref{it2}) imply

$$p(a,b+1,c)-p(a,b,c)=-{\frac{1}{D(a,c)}}U_bU_{b-1}\cdots U_{a+1}
.$$

\noindent Adding these equations we obtain

\begin{eqnarray*}
p(a,b+1,c)-1&=&p(a,b+1,c)-p(a,a,c)=\\
&=&-{\frac{1}{D(a,c)}}(1+U_{a+1}+U_{a+1}U_{a+2}+\cdots+
U_{a+1}U_{a+2}\cdots U_b)=\\
&=&-{\frac{D(a,b+1)}{D(a,c)}}.
\end{eqnarray*}

\noindent Hence we have the lemma. $\Box$

Introduce the following notations:

\begin{eqnarray*}
\lambda(0,i)&=&1,\\
\lambda(1,i)&=&i,\\
\lambda(2,i)&=&\lambda(1,i)\log i,\ldots,\\
\lambda(k,i)&=&\lambda(k-1,i)\log_{k-1}i\quad (k=3,4,\ldots),
\end{eqnarray*}

\noindent where

\begin{eqnarray*}
\log_0i&=&i,\\
\log_1i&=&\log i,\ldots,\\
\log_ki&=&\log\log_{k-1}i,
\end{eqnarray*}

\noindent and

\begin{eqnarray*}
\Lambda(0,i)&=&0,\\
\Lambda(K,i)&=&\sum_{k=1}^K{\frac{1}{\lambda(k,i)}},\quad (K=1,2,
\ldots)\\
\Lambda(K,i,B)&=&\Lambda(K-1,i)+{\frac{B}{\lambda(K,i)}} \quad
(B>0).
\end{eqnarray*}

\noindent Note that

\begin{eqnarray*}
\Lambda(1,i,B)&=&{\frac{B}{i}},\\
\Lambda(2,i,B)&=&{\frac{1}{i}}+{\frac{B}{i\log i}},\\
\Lambda(3,i,B)&=&{\frac{1}{i}}+{\frac{1}{i\log i}}+{\frac{B}
{i\log i\log\log i}}.
\end{eqnarray*}

\noindent For some $K=1,2,\ldots,\ B>0$   f\/ixed, define

$$i_0=\min \left\{i:\,\frac{1}{4}\Lambda(K,i,B)<\frac{1}{2}\right\}$$
and let
\begin{equation}
p_i=\left\{\begin{array}{ll}
p_{i_0},\quad & i\!f \quad 1\leq i\leq i_0\\
 \frac{1}{4}\Lambda(K,i,B) \quad & i\!f \quad i>i_0.
\end{array}\right. \label{pei}
\end{equation}
 \noindent Now
we are interested in the case $\{p_i\}$ above. In fact, in the
future  for convenience, when we say that
$$\displaystyle {p_i=\frac{1}{4}\Lambda(K,i,B)}$$
we actually mean that $p_i$ is defined by (\ref{pei}).
\begin{lemma}

\noindent Let $\displaystyle {p_i=\frac{1}{4}\Lambda(K,i,B)}$.
Then

\begin{equation}
D(0,\infty)\left\{\begin{array}{ll}
\label{conv}=\infty\quad & i\!f \quad B\leq 1,\\
<\infty\quad & i\!f\quad B>1,
\end{array}\right.
\end{equation}

\begin{equation}
p(0,1,\infty)\left\{\begin{array}{ll}
\label{conv2}=1\quad & i\!f \quad B\leq 1,\\
<1\quad & i\!f\quad B>1.
\end{array}\right.
\end{equation}
\noindent For $n\geq m+2,$  $B\neq1$ and $m$ big enough
\begin{eqnarray}
D(m,n)&=&(1+o_m(1))\lambda(K-1,m)(\log_{K-1}m)^B\times \label{dmn}\\
\nonumber
&\times&\frac{1}{B-1}\left({\frac{1}{(\log_{K-1}m)^{B-1}}}
-{\frac{1}{(\log_{K-1}n)^{B-1}}}\right).
\end{eqnarray}

\noindent  If $B>1$,

\begin{equation}
D(m,\infty)=(1+o_m(1)){\frac{\lambda(K,m)} {B-1}}, \label{demi}
\end{equation}

\begin{equation} \label{mmpl}
p(m,m+1,\infty)=1-(1+o_m(1)){\frac{(B-1)}{\lambda(K,m)}}.
\end{equation}
\end{lemma}

\noindent
{\bf Proof}:  To prove (\ref{dmn}), observe that  from (\ref{ui})
we have for $n\geq m+2$

\begin{eqnarray}D(m,n)&=& \label{dmn2}
1+\sum_{j=1}^{n-m-1}\prod_{i=1}^j U_{m+i}\\&= &\nonumber
1+\sum_{j=1}^{n-m-1} \exp\left(-\sum_{i=m+1}^{m+j}\left(
\Lambda(K,i,B)\right)\right) \exp\left(O(1)
\sum_{i=m+1}^{m+j}\Lambda^2(K,i,B)\right)\\&=&\nonumber
1+(1+o_m(1))\sum_{j=1}^{n-m-1}\exp\left(-\sum_{i=m+1}^{m+j}
\Lambda(K,i,B)\right)\\&=:&\nonumber 1+(1+o_m(1))A(m,n,K).
\end{eqnarray}
Now we give a lower bound for $A(m,n,K).$

\begin{eqnarray}
A(m,n,K)&\geq&\label{dmn3}
\sum_{j=1}^{n-m-1}\exp\left(-\int_{m}^{m+j}\Lambda(K,x,B)\,dx\right)\\
&=&\nonumber \sum_{j=1}^{n-m-1}\frac{\lambda(K-1,m)
(\log_{K-1}m)^B}{\lambda(K-1,m+j) (\log_{K-1}(m+j))^B}\\
 &=&\nonumber \lambda(K-1,m)
(\log_{K-1}m)^B\sum_{\ell=m+1}^{n-1}\frac{1}
{\lambda(K-1,\ell)(\log_{K-1}\ell)^B}
\\&\geq&\nonumber \lambda(K-1,m)
(\log_{K-1}m)^B\int_{m+1}^{n}\frac{1}{\lambda(K-1,x)(\log_{K-1}x)^B}\,dx
\\&=& \nonumber\lambda(K-1,m)
(\log_{K-1}m)^B\left(\frac{(\log_{K-1}m)^{1-B}-
(\log_{K-1}n)^{1-B}}{B-1}\right).
\end{eqnarray}
It is easy to see that the proof of the upper bound goes the same way,
resulting the same expression as in (\ref{dmn3}) with $m$
replaced by $m+1$ which combined with (\ref{dmn2}) proves (\ref{dmn}).
The proof of (\ref{conv}) is similar, and the rest of the lemma
follows from these two. $\Box$
\newpage
\noindent
 {\bf Consequence}: If  for any $K=1,2...$

$$\displaystyle {p_i=\frac{\Lambda(K,i,B)}{4}},$$
then the Markov chain is recurrent if $B\leq 1$ and transient if
$B>1. $

Now we would like to consider the case when $p_i$ is essentially
$\displaystyle{\frac{B}{4i^{\alpha}}}$, which should be understood
in the same way as it was defined in (\ref{pei}). Namely,  let
$$i_0=\min \left\{i:\,\frac{B}{4i^{\alpha}}<\frac{1}{2}\right\}$$
and let
\begin{equation}
p_i=\left\{\begin{array}{ll}
 \displaystyle{p_{i_0}},\quad & i\!f \quad 1\leq i\leq i_0\\
 \displaystyle{\frac{B}{4i^{\alpha}}} \quad & i\!f \quad i>i_0.
\end{array}\right.
\end{equation}

\vspace{2ex}\noindent
\begin{lemma}
In case $\displaystyle{p_i=\frac{B}{4i^{\alpha}}\ (0<\alpha<1)}$
we have

\begin{equation} \label{dal}
D(m,\infty)=(1+o_m(1))\frac{m^\alpha}{B},
\end{equation}
\begin{equation}\label{pmm}
1-p(m,m+1,\infty)=(1+o_m(1))\frac{B}{m^{\alpha}}.
\end{equation}
\end{lemma}

  \vspace{2ex}\noindent
{\bf Proof:}  Consider the case $0<\alpha<1/2$ first.  By
  (\ref{ui})
  \begin{eqnarray*}
 && \prod_{i=1}^j U_{m+i}\leq\exp\left(-B\sum_{\nu=m+1}^{m+j}
   \nu^{-\alpha}+\sum_{\nu=m+1}^{m+j}
  C\nu^{-2\alpha}\right)=\\
  &\leq&(1+o_m(1))\exp\left(\frac{-B}{1-\alpha}\left[(m+j)^{1-\alpha}-
  (m)^{1-\alpha}\right]+\frac{C}{1-2\alpha}\left[(m+j)^{1-2\alpha}-
  (m)^{1-2\alpha}\right]\right).
  \end{eqnarray*}
 \noindent  Consequently,
\begin{eqnarray*}
&&D(m,n)\\&\leq&(1+o_m(1))\exp\left(\frac{B
m^{1-\alpha}}{1-\alpha}- \frac{Cm^{1-2\alpha}}{1-2\alpha}
\right)\sum_{k=m+1}^n
\exp(\left(-\frac{Bk^{1-\alpha}}{1-\alpha}+\frac{C
k^{1-2\alpha}}{1-2\alpha}
\right)\\&\leq&(1+o_m(1))\exp\left(\frac{B
m^{1-\alpha}}{1-\alpha}- \frac{Cm^{1-2\alpha}}{1-2\alpha}\right)
\int_{m+1}^n
\exp\left(-\frac{Bx^{1-\alpha}}{1-\alpha}(1-C\frac{1-\alpha}{1-2\alpha}
x^{-\alpha})\right)
\,dx\\&\leq&(1+o_m(1))\exp\left(B\frac{m^{1-\alpha}}{1-\alpha}-
\frac{Cm^{1-2\alpha}}{1-2\alpha}\right) \int_{m+1}^n
\exp\left(-\frac{Bx^{1-\alpha}}{1-\alpha}h_m\right)\,dx,
\end{eqnarray*}
 \noindent  where  $$ h_m=1-C\frac{1-\alpha}{1-2\alpha}(m+1)^{-\alpha}.$$
In the calculation above $C$ is a positive constant the value of
which is not important. In the future we will use $C,\,C^*$ or
$C_1, C_2 \ldots$ for which this remark applies, and their values
might change from line to line.  Using substitution and the
asymptotic representation of the incomplete Gamma function (see
e.g. Gradsteyn and Ryzhik \cite{RG} page 942, formula (8.357))
$$\Gamma(\beta,x)=\int_x^{\infty}t^{\beta-1}e^{-t}\,dt=x^{\beta-1}e^{-x}\left(
1+\frac{O(1)}{x}\right)\qquad {\rm as}\,\, x\to \infty$$ we
conclude that as $m\to \infty$

\begin{eqnarray*}
D(m,\infty)&\leq&\left(1+O(\frac{1}{m^{1-\alpha}})\right)\frac{m^\alpha}{Bh_m}
\exp\left(\frac{Bm^{1-\alpha}}{1-\alpha}-\frac{Cm^{1-2\alpha}}{1-2\alpha}
\right) \exp\left(-\frac{h_m B m^{1-\alpha}}{1-\alpha}\right)
\\&=&\left(1+O(\frac{1}{m^{1-\alpha}})\right)\frac{m^{\alpha}}{B}.
\end{eqnarray*}
A similar calculation (which we omit) gives the same lower bound.
The case of $\alpha=1/2 $ goes along the same lines with obvious
modifications. On the other hand, the case $1/2<\alpha<1$ can be
worked out similarly, but it   is obvious with  less precise
calculations  as well. $\Box$
\begin{lemma}
In case $\displaystyle{p_i=\frac{B}{4 (\log i)^{\al}}}$ with
$\alpha>0,$ there exist $0<K_1<K_2$ such that
\begin{equation}
K_1(\log m)^{\al}\leq D(m,\infty)\leq K_2(\log m)^{\al},
\end{equation}
\begin{equation}\label{log}
1-p(m,m+1,\infty)=\frac{O(1)}{(\log m)^{\al}}.
\end{equation}
\end{lemma}

\noindent
{\bf Proof:} First we give the upper bound. For $m\geq m_0$
\begin{eqnarray*}
\sum_{i=m}^{m+j}\left(\frac{B}{(\log i)^{\al}}-\frac{C}{(\log
i)^{2\al}}\right)&=&\sum_{i=m}^{m+j}\frac{B}{(\log
i)^{\al}}\left(1-\frac{C^*} {(\log i)^{\al}}\right)\\ &\geq&
\sum_{i=m}^{m+j}\frac{B(1-\ep)}{(\log
i)^{\al}}=:A(m,j,\ep).\end{eqnarray*}

\noindent Then for
 $$\ell(\log m)^{\al}\leq j<(\ell+1)(\log m)^{\al}\quad (\ell=0,1,2,\ldots)$$
we have
$$
 A(m,j,\ep)\geq{\frac{B(1-\ep)\ell(\log m)^{\al}}{(\log[m+(\ell+1)
(\log m)^{\al}])^{\al}}}=:H(m,\ell,\alpha) .
$$
It is easy to see now, that if $(\ell+1)(\log m)^{\al} \leq m$
then for an appropriate $C_1$
$$H(m,\ell,\alpha)\geq  \frac{B(1-\ep)\ell(\log m)^{\al}}{(\log (2m))^\alpha}
 \geq C_1 \ell.$$
On the other hand, if $(\ell+1)(\log m)^{\al} \geq m,$ then for an
appropriate $C_2$
$$H(m,\ell,\alpha)\geq  \frac{B(1-\ep)\ell(\log m)^{\al}}
{(\log (2(\ell +1)(\log m)^{\alpha})^{\alpha}}
 \geq C_2 {\ell}^{1/(2\alpha)}.$$
 \noindent Then  with $N=N(\alpha):=[\frac{m}{(\log m)^{\alpha}}]-1.$

\begin{eqnarray*}
D(m,\infty)&\leq&\sum_{\ell=0 }^N   e^{-C_1\ell}\,(\log
m)^{\al}+\sum_{\ell=N }^{\infty}
e^{-C_2\displaystyle{\ell^{\frac{1}{2\alpha}}}}(\log
m)^{\al}=O(1)(\log m)^{\al}.
\end{eqnarray*}
The lower bound follows from Lemma 2.1. $\Box$

\section{Local time}

\setcounter{equation}{0}
 \setcounter{lemma}{0}
\setcounter{theorem}{0}

\noindent We intend to study the limit properties of the local
time $\xi(R,\infty)$ in case of transient random walks. To this
end we also define the  number of upcrossings by

\begg \xi(R,n,\uparrow):=\#\{k:\ 0\leq k\leq n,\ X_k=R,\
X_{k+1}=R+1\}.\label{upcr}\endd

\begg \xi(R,\infty,\uparrow):=\lim_{n\to
\infty}\xi(R,n,\uparrow).\endd

\begin{lemma} For $R=0,1,2,\ldots $
\begin{equation}  \label{loc}
{\bf P}(\xi(R,\infty)=L)={\frac{1+2p_R}{2D(R,\infty)}}
\left(1-{\frac{1+2p_R}{2D(R,\infty)}}\right)^{L-1},\quad L=1,2,\ldots
\end{equation}
Moreover, the sequence
$$
\xi(R,\infty,\uparrow),\qquad R=0,1,2,\ldots
$$
is a Markov chain and
\begin{equation}
{\bf P}(\xi(R,\infty,\uparrow)=L)=\frac{1}{D(R,\infty)}
\left(1-\frac{1}{D(R,\infty)}\right)^{L-1},\quad L=1,2,\ldots
\label{up}
\end{equation}
\end{lemma}

\noindent {\bf Proof:} Clearly we have for $L=1,2, \ldots$
\begin{eqnarray*}
\lefteqn{\bp(\xi(R,\infty)=L)=
\left({\frac{1}{2}}+p_R\right)(1-p(R,R+1,\infty))
\times}\\
&&\times\sum_{j=0}^{L-1}{{L-1}\choose{j}}\left({\frac{1}{2}}-p_R\right)^j
\left(\left({\frac{1}{2}}+p_R\right)p(R,R+1,\infty)\right)^{L-j-1}=\\
&=&\left({\frac{1}{2}}+p_R\right)(1-p(R,R+1,\infty))\times\\
& &\times\left(1-\left({\frac{1}{2}}+
p_R\right)(1-p(R,R+1,\infty))\right)^{L-1},
\end{eqnarray*}
implying (\ref{loc}) by (\ref{it0}).

 The other statements of the Lemma are obvious. $\Box$

\begin{theorem}If $p_R\to 0$, as $R\to\infty$, then

$$
\lim_{R\to \infty} \bp \left(
\frac{\xi(R,\infty)}{ 2 D(R,\infty)}>x \right)=
\lim_{R\to \infty} \bp \left(
\frac{\xi(R,\infty,\uparrow)}{D(R,\infty)}>x \right)=
e^{-x},
$$
that is to say, $\displaystyle{ \frac{\xi(R,\infty)}{ 2
D(R,\infty)}}$ and $\displaystyle{ \frac{\xi(R,\infty,\uparrow)}
{D(R,\infty)}}$
have exponential limiting distributions.
\end{theorem}

\noindent The proof is a trivial consequence of Lemma 3.1.
 \vspace{2ex}\noindent

\begin{theorem}

Assume that $p_R\to 0$ as $R\to\infty$. Then with probability $1$ we
have

\begg
  \xi(R,\infty)\leq 2(1+\varepsilon)
D(R,\infty)\log R
\label{loc2}\endd
for any $\ep>0$ if $R$ is large enough.

Moreover,
 \begg
 \xi(R,\infty)\geq M D(R,\infty) \quad  { \rm i.o.\,\,a.s.}
{\label{loc3}}\endd for any $M>0.$
\end{theorem}
In case $\displaystyle{p_R=\frac{\Lambda(K,R,B)}{4}}$  with
$B>1,$ instead of  (\ref{loc2}) and (\ref{loc3}) we have the much
sharper

\begin{theorem}
For $\displaystyle{p_R=\frac{\Lambda(K,R,B)}{4}}$, $B>1$, we have
\begg
\label{lilone} \limsup_{R\to \infty} \frac{\xi(R,\infty)}{2
D(R,\infty)\log\log R}\leq 1.
\endd
and
\begg \label{lilone*} \limsup_{R\to \infty}
\frac{\xi(R,\infty)}{2 D(R,\infty)\log_{K+1} R}\geq 1.
\label{Kone}
\endd \noindent
 Especially  in case
$\displaystyle{p_R=\frac{\Lambda(1,R,B)}{4}}=\frac{B}{4R},$ $B>1$,
being  $\displaystyle{D(R,\infty)=\frac {R}{B-1}}$, we have
 \begin{equation} \label{liltwo}
 \limsup_{R\to \infty} \frac{(B-1)\xi(R,\infty)}{2
R\log \log R}=\limsup_{R\to \infty} \frac{(B-1)\xi(R,\infty,\uparrow)}
{R\log \log R}= 1. \end{equation}
\end{theorem}

\noindent {\bf Consequences:}
\begin{itemize}
\item If $\displaystyle{p_R=\frac{1}{4}\Lambda(K,R,B),\ (B>1)}$
then for any $\varepsilon> 0 $ \noindent

\begin{equation} \label{loc11}
\xi(R,\infty)\leq
{\frac{2(1+\varepsilon)}{B-1}}\lambda(K,R)\log\log R {\qquad \rm
a.s.}
\end{equation}
if $R$ is large enough
\begin{equation} \label{loc111}
\xi(R,\infty)\geq{\frac{2(1-\varepsilon)}{B-1}}\lambda(K,R)\log_{K+1}R
{\qquad \rm i.o.\,\,a.s.}\end{equation}

\noindent and
\begin{equation}
\lim_{R\to \infty} \bp \left( \frac{B-1}{ 2
\lambda(K,R)}\,\xi(R,\infty)>x \right)=e^{-x}.
\end{equation}
 \item If $\displaystyle{p_R=\frac{B}{4 R^{\alpha}}\ (0<\alpha<1),}$
\noindent then

\begin{eqnarray}\label{al1}
\xi(R,\infty)&\leq&\frac{2}{B}(1+\varepsilon)R^\alpha\log R\quad
{\rm a.s.},\\
\label{al2} \xi(R,\infty)&\geq& MR^\alpha\quad {\rm i.o.\,\,a.s.}
\end{eqnarray}
for any $M>0$ and

\begin{equation}
\lim_{R\to \infty} \bp \left( \frac{B\,
\xi(R,\infty)}{2R^{\alpha}}>x \right)=e^{-x}.
\end{equation}

\item If $\displaystyle{p_R=\frac{B}{4(\log R)^{\alpha}}\quad
(\alpha>0),}$ then

\begin{eqnarray}\label{logal}
\xi(R,\infty)&\leq&O(1)(\log R)^{1+\alpha}\quad {\rm a.s.},\\
\label{logal2} \xi(R,\infty)&\geq& M(\log R)^{\alpha}\quad {\rm
i.o.\,\,a.s.}
\end{eqnarray}
for any $M>0$.
\end{itemize}

\vspace{2ex}\noindent \noindent {\bf Proof of Theorem 3.2:}
(\ref{loc2}) follows from Lemma 3.1. \noindent On the other hand,
(\ref{loc}) also implies that for any $M>0$

$$\liminf_{R\rightarrow\infty}\bp\left(\xi(R,\infty)\geq M
D(R,\infty)\right)>0.$$

Now to finish our proof  we need to apply the zero-one law (in a
non-independent setup) exactly in the same way as in Doob
\cite{DO} page 103, observing   that the conditional probability
of our tail event given the first $n$ steps of our walk is the
same as its unconditional probability, that is for any $n=1,2,\ldots$
$$\bp\left(\xi(R,\infty)\geq M
D(R,\infty)\,\,{\rm i.o.}\mid X_1,X_2, ...X_n\right)=
\bp\left(\xi(R,\infty)\geq M D(R,\infty)\,\,{\rm i.o.}\right).
$$
\noindent which, in turn, implies (\ref{loc3}).

\vspace{2ex}\noindent\noindent {\bf Proof of Theorem 3.3:}

To prove  (\ref{lilone}), we need a few lemmas. Recall the
definition of the upcrossing in (\ref{upcr}). For large values of
the local time and upcrossing we have the following invariance
principle.

\begin{lemma} As $R\to\infty$
\begin{equation}
\xi(R,\infty)-2\xi(R,\infty,\uparrow)=O((D(R,\infty) \log R
)^{1/2+\varepsilon} +p_R D(R,\infty)\log R)\quad {\rm a.s.}
\label{inv}
\end{equation}
\end{lemma}

\noindent{\bf Proof:} Under the condition $\xi(R,\infty)=L$,
$\xi(R,\infty,\uparrow)-1$ has binomial distribution with
parameters $(L-1,1/2+p_R)$. According to Hoeffding inequality,
$$
\bp\left(\left|\xi(R,\infty,\uparrow)-1-
\left(\frac12+p_R\right)(L-1)\right|\geq u(L-1)^{1/2}\right)\leq
e^{-Cu^2}
$$
with some $C>0$, from which as $L\to\infty$,
$$
\xi(R,\infty,\uparrow)-\frac{L}{2}=O(L^{1/2+\varepsilon}+Lp_R)\quad{\rm
a.s.}
$$
Putting $L=\xi(R,\infty)$,  we get (\ref{inv}) from (\ref{loc2}).

\begin{lemma} Let
$$
\gamma_R=\left(\frac12 +p_R\right)p(R,R+1,\infty),
$$
and
$$c_R=\frac{\gamma_R}{1-\gamma_R}.$$
Then
$$
\zeta(R):=\frac{\xi(R,\infty,\uparrow)}{c_1\cdots c_R},\quad
R=1,2,\ldots
$$
is a submartingale.
\end{lemma}

\noindent{\bf Proof:}
Let $T_R$ be the first hitting time of $R$ by $\{X_n\},$
e.g. $T_R=\min\{n: X_n=R\}.$ Then we have

\begin{equation}\bp_R(\xi(R,T_{R-1},\uparrow)=j, T_{R-1}<\infty)=
\left(\frac12
-p_R\right)\gamma_R^j,\quad j=0,1,\ldots, \label{upr1}
\end{equation}

\begin{equation}
\bp_R(\xi(R,\infty,\uparrow)=j,T_{R-1}=\infty)= \left(\frac12
+p_R-\gamma_R\right)\gamma_R^{j-1},\quad j=1,2,\ldots \label{upr2}
\end{equation}
Observe that
$$
\xi(R,\infty, \uparrow)=\sum_{m=1}^{\xi(R-1,\infty,
\uparrow)-1}\xi_m +\widetilde\xi,
$$
where $\xi_m, \,\, m=1,2...$ has distribution (\ref{upr1}) and
$\widetilde\xi$ has distribution (\ref{upr2}). Then \begg
\be(e^{\lambda\xi(R,\infty, \uparrow)},\,
\xi(R-1,\infty,\uparrow)=i)
=(\be(e^{\lambda\xi_1}))^{i-1}\be(e^{\lambda\widetilde\xi}) \label
{gener}\endd
$$
=\frac{\left(\frac12 +p_R-\gamma_R\right)e^\lambda
\left(\frac12-p_R\right)^{i-1}}{(1-\gamma_R e^\lambda)^i},
$$
hence
$$
\be (e^{\lambda\xi(R,\infty, \uparrow)}\mid
\xi(R-1,\infty,\uparrow)=i)
=e^\lambda\left(\frac{1-\gamma_R}{1-\gamma_R e^\lambda}\right)^i,
$$
from which
\begg \be (\xi(R,\infty,\uparrow)\mid \xi(R-1,
\infty,\uparrow))=c_R\xi(R-1,\infty,\uparrow)+1,
 \label{upr3}\endd
which easily implies the lemma. $\Box$

Now we prove the upper bound, i.e.
\begin{equation}
\limsup_{R\to\infty}\frac{\xi(R,\infty,\uparrow)}{D(R,\infty)\log\log
R}\leq 1 \qquad {\rm a.s.}, \label{lil}
\end{equation}
which  also implies  (\ref{lilone}) by Lemma 3.2.

With an easy calculation we get  from (\ref{gener}) that \begg
\be(e^{\lambda\xi(R,\infty,\uparrow)})=
\frac{e^\lambda}{D(R,\infty)-e^{\lambda}(D(R,\infty)-1)}.
\label{mgf}
\endd

Using that $\zeta(R)$ is submartingale, from (\ref{mgf}) we have with
$R_k=[\exp(k/\log k)]$, $C_k=c_1c_2\ldots c_{R_k}$,
$$
u_k=(1+\varepsilon)D(R_k,\infty)\log\log R_k,
$$

$$
\bp\left(\max_{R_k\leq R< R_{k+1}}\zeta (R)\geq
\frac{u_{k+1}}{C_{k+1}}\right)
$$
$$
\leq \exp(-\lambda u_{k+1}/C_{k+1})\be(\exp(\lambda\zeta(R_{k+1})))
$$
$$
=\frac{\exp(\lambda/C_{k+1})(1-u_{k+1})}
{D(R_{k+1},\infty)-\exp(\lambda/C_{k+1})(D(R_{k+1},\infty)-1)}.
$$

It can be seen that the optimal choice for $\lambda$ is given by
$$
\exp(\lambda/C_{k+1})=\frac{(u_{k+1}-1)D(R_{k+1},\infty)}
{u_{k+1}(D(R_{k+1},\infty)-1)},
$$
and we get finally
$$
\bp\left(\max_{R_k\leq R< R_{k+1}}\zeta (R)\geq
\frac{u_{k+1}}{C_{k+1}}\right)
=\frac{O(1)\log\log R_{k+1}}{(\log R_{k+1})^{1+\varepsilon}}.
$$

Hence by Borel-Cantelli lemma for large $k$ and $R_k\leq R<R_{k+1}$
we have
$$
\zeta(R)\leq \frac{(1+\varepsilon)D(R,\infty)\log\log R}
{c_1\cdots c_R c_{R+1}\cdots c_{R_{k+1}}},
$$
i.e.
$$
\xi(R,\infty,\uparrow)\leq \frac{(1+\varepsilon)
D(R,\infty)\log\log R} {c_{R+1}\cdots c_{R_{k+1}}}.
$$

If $p_R=\Lambda(K,R,B)/4$, then (cf. (\ref{demi}))
$$
D(R,\infty)\sim\frac{\lambda(K,R)}{B-1}
$$
and
$$
c_R\sim \frac{1+2p_R-1/D(R,\infty)}{1-2p_R+1/D(R,\infty)}
\sim \exp(4p_R-2/D(R,\infty))\sim\exp\left(\Lambda(K,R,B)
-\frac{2(B-1)}{\lambda(K,R)}\right).
$$

If $K=1$, then
$$
\Lambda(1,R,B)-\frac{2(B-1)}{\lambda(1,R)}\sim \frac{2-B}R, \quad B\neq
2
$$
and
$$
\Lambda(1,R,2)-\frac{2}{\lambda(1,R)}=\frac{o(1)}R,
$$
and if $K>1$, then
$$
\Lambda(K,R,B)-\frac{2(B-1)}{\lambda(K,R)}\sim\frac1R.
$$
Hence for large $k$ and $R_k\leq R\leq R_{k+1}$ we have
$$
c_{R+1}\cdots c_{R_{k+1}}\sim \exp\left(C\log \frac{R_{k+1}}R\right)
$$
with some constant $C$ if $K=1,\, B\neq 2$ or $K>1$ and $C=o(1)$ if
$K=1, \, B=2$. In view of $\lim_{k\to\infty}R_{k+1}/R_k=1$,
for any $\varepsilon>0$, one can choose $k$ large enough such that
$$
c_{R+1}\cdots c_{R_{k+1}}\geq 1-\varepsilon,
$$
i.e.
$$
\xi(R,\infty, \uparrow)\leq \frac{(1+\varepsilon)
D(R,\infty)\log\log R} {1-\varepsilon}.
$$
Since $\varepsilon>0$ is arbitrary, (\ref{lil}) follows.

To prove the lower bound (\ref{Kone}), consider an
increasing sequence of sites $R_k$ to be determined later.
 Let
$$\tau_k=\min\{n: X_n=R_k\},$$
the time of the first visit at the site $R_k,$ and define
$$Z(k):=\xi(R_k,\tau_{k+1}).$$
\noindent Observe that  $ \{Z(k),  \, k=1,2...\}$ are independent.
Following the proof of Lemma 3.1 we can conclude that
\begin{eqnarray}
\lefteqn{\bp(Z(k)\geq L)=(1+o_{R_k}(1))\times} \nonumber \\
&\times&\left[\left(1-{\frac{1}{2}}(1-
p(R_k,R_k+1,R_{k+1}))\right)\left(
1+O((1-p(R_k,R_k+1,R_{k+1}))p_{R_k})\right)\right]^{L-1}.
\end{eqnarray}
\noindent  Based on (\ref{dmn}) it is easy to calculate that
\begin{eqnarray}
  D(R_k,R_{k+1})&=& \left(1+o_{R_k}(1)\right)
\frac{\lambda(K-1,R_k)}{B-1}\log_{K-1}R_k\left(1-\left(
\frac{\log_{K-1}R_k}{\log_{K-1}R_{k+1}}\right)\right)\nonumber=\\
&=&\left(1+o_{R_k}(1)\right)\frac{\lambda(K,R_k)}{B-1}\left(1-\left(
\frac{\log_{K-1}R_k}{\log_{K-1}R_{k+1}}\right)\right).\end{eqnarray}
Define the sequence $R_k$ by
 $$ \log_K R_k:=k\log Q$$
  with some $Q>1$  (we intentionally  forget about the technicalities
arising from the fact that
  the sites should be integers).
  It is easy to see that with this choice of $R_k$
$$\frac{\log_{K-1}R_k}{\log_{K-1}R_{k+1}}=\frac{1}{Q}.$$

Let
$$L(k)=2\frac{\lambda( K,R_k)}{B-1}\frac{Q-1}{Q}\log_{K+1} R_k.  $$
 From (\ref{it2}) we get that
$$\bp (Z(k)\geq L(k))\sim \exp(-\log_{K+1}R_k)=\frac{1}{\log_K R_k}
=\frac{1}{k\log Q}.$$

   \noindent Applying Borel-Cantelli lemma and then
letting $Q\to \infty$, we get (\ref{Kone}). $\Box$

Our next issue was to investigate how small could be the  local
time of our process. More precisely we wanted  to know whether it
is true that in the transient case there are always infinitely
many sites with local time equal to $1$. In fact we managed to
prove in some sense much more, and in some sense much less.
Namely, we  prove the following two theorems. Define for $N\geq 2$

$$f(N,R)=f(N,R,\varepsilon)=\frac{1}{\log 2}
\left(\sum_{j=2}^{N}\log_{j} R+{\varepsilon \log_{N}R}\right)$$
and $$g(N,R)=f(N,R,0).$$
\begin{theorem} Let $\displaystyle{p_R=\frac{\Lambda(1,R,B)}{4}}$ with $ B>1$ and
$N\geq 2$. Then
\begin{itemize}
\item
  with probability $1$ there exist infinitely many $R$ for which

$$\xi(R+j,\infty)=1$$

\noindent for each $j=0,1,2,\ldots,[g(N,R)]$. \item

   with probability $1$ for any $\varepsilon>0$ and $R$ large enough
 there exists an $S$

$$R\leq S\leq f(N,R,\ep)$$

such that

$$\xi(S,\infty)>1.$$
\end{itemize}
\end{theorem}

Let $$f^*(R,\epsilon)=\frac{(1+\varepsilon)(1-\alpha)\log R}{\log
2}\qquad {\rm and}\qquad g^*(R)=f^*(R,0)$$

 \noindent
\begin{theorem} Let $\displaystyle{p_R=\frac{B}{4 R^{\alpha}}\ (0<\alpha<1)}$.
 Then
\begin{itemize}
\item
 with probability $1$ there exists infinitely many $R$ for which

$$\xi(R+j,\infty)=1$$

for each $j=0,1,2,\ldots,g^*(R) $. \item  with probability $1$
for each $R$ large enough and $\varepsilon>0$ there exists an $S$,

$$R\leq S\leq f^*(R,\varepsilon)$$

such that

$$\xi(S,\infty)>1.$$
\end{itemize}
\end{theorem}

Furthermore, we conjecture that for $p_i \geq B/(4i)$, where $B>1,$
with probability $1$ there are always infinitely many sites with
local time $1.$ On the other hand, recently James et al.
\cite{JLP} proved that for $p_i\sim\Lambda(2,i,B) $ with $B>1$
with probability $1$ there are only finitely many cutpoints, hence
finitely many points with local time $1.$ We note that it can be
seen with  a similar argument that this is the case for $p_i\sim
\Lambda(K,i,B) $ for all $K\geq 2$ as well.

 {\bf Proof of Theorem 3.4:} At f\/irst we prove the second
statement. Recall the notation of $\lambda(N,R)$ and observe
that \begg R2^{g(N,R)}=\lambda(N,R)\quad {\rm and}\quad
R2^{f(N,R)}=\lambda(N-1,R)
(\log_{N-1}R)^{1+\epsilon}.\label{nohat}\endd

 Now the proof of the second statement is a trivial consequence of

\begin{lemma} For every $N\geq 2$ integer as $R\to \infty$

\begin{eqnarray*}
\lefteqn{\bp\left\{\bigcap_{j=1}^{f(N,R)}\{\xi(R+j,\infty)=1\}\right\}=}\\
&=&\prod_{j=1}^{f(N,R)}\left({\frac{1}{2}}+{\frac{B}{4(R+j)}}\right)
(1-p(R+f(N,R),R+f(N,R)+1,\infty))=\\
&=&(1+o_R(1)){\frac{1}{2^{f(
N,R)}}}{\frac{B-1}{R}}=(1+o_R(1)){\frac{B-1} {\lambda
(N-1,R)(\log_{N-1} R)^{1+\varepsilon}}}.
\end{eqnarray*}
\end{lemma}

 \vspace{2ex}\noindent
{\bf Proof:} Obvious by (\ref{mmpl}). $\Box$

The proof of the first statement of the theorem is based on the
following

\begin{lemma} For every $N\geq 2$ integer as $R\to \infty$

\begg
\bp\left\{\bigcap_{j=1}^{g(N,R)}\{\xi(R+j,\infty)=1\}\right\}=
{\frac{O(1)}{\lambda(N,R)}},\label{g1} \endd

\begin{eqnarray}\label{g2}
\lefteqn{{\cal P}:={\cal P}(N,R,S):=}\\
&=&\bp\left\{\bigcap_{j=1}^{g(N,R)}\{\xi(R+j,\infty)=1\}\cap\bigcap_{j=1}
^{g(N,S)}\{\xi(S+j,\infty)=1\}\right\}\leq\nonumber\\
&\leq&\left\{\begin{array}{ll}
\displaystyle{{\frac{(1+o_R(1))(B-1)^2}{\lambda(N,R)\lambda(N,S-R)}}}\
& i\!f\ S\geq
R+g(N,R),\vspace{2ex}\\
\displaystyle{{\frac{O(1)2^R}{2^{S+g(S,N)}}}{\frac{B-1}{S+g(N,S)}}}\
& i\!f \ R<S<R+g(N,R).
\end{array}\right.\nonumber
\end{eqnarray}
\end{lemma}

\vspace{2ex}\noindent {\bf Proof:} (\ref{g1}) follows from Lemma
3.1 and  (\ref{nohat}). In case $R<S<R+g(N,R)$ we have

\begin{eqnarray*}
{\cal
P}&=&\prod_{i=R}^{S+g(N,S)}\left({\frac{1}{2}}+{\frac{B}{4i}}\right)
(1-p(S+g(N,S),S+g(N,S)+1,\infty))\leq\\
&\leq& O(1){\frac{1}{2^{S+g(N,S)-R}}}{\frac{B-1}{S+g(N,S)}}.
\end{eqnarray*}

 \noindent In case $S>R+g(N,R)$ we have

\begin{eqnarray*}
{\cal P}&=&(1+o_R(1)){\frac{1}{2^{g(N,R)}}}\,(1-p(R,R+1,S))
{\frac{1}{2^{g(N,S)}}}(1-p(S,S+1,\infty))=\\
&=&(1+o_R(1)){\frac{1}{2^{g(N,R)}\,2^{g(N,S)}}}{\frac{B-1}{R^B}}{\frac{1}
{R^{1-B}-S^{1-B}}}{\frac{B-1}{S}}=\\
&=&(1+o_R(1)){\frac{1}{2^{g(N,R)}\,2^{g(N,S)}}}
{\frac{(B-1)^2}{R}}{\frac{S^{B-2}}{S^{B-1}-R^{B-1}}}\leq\\
&\leq&(1+o_R(1)){\frac{1}{2^{g(N,R)}\,2^{g(N,S-R)}}}
{\frac{(B-1)^2}{R}}{\frac{1}{(S-R)}} \leq \\
&\leq&(1+o_R(1)){\frac{(B-1)^2}{\lambda(N,R)\lambda(N,S-R)}}.
\end{eqnarray*}

\noindent Hence we have the second statement of the lemma. $\Box$

\vspace{2ex}\noindent Now we turn to the proof of the first
statement  of the theorem. Let

$$A(R)=\bigcap_{j=1}^{g(N,R)}\{\xi(R+j,\infty)=1\}.$$

\noindent Then by (\ref{g1})

\begin{equation}\label{g3}
\sum_{R=1}^T\bp(A(R))=O(1)\log_N T
\end{equation}

\noindent and

\begin{eqnarray*}
\lefteqn{\sum_{R=1}^T\sum_{S=R+1}^T\bp(A(R)A(S))=}\\
&=&\sum_{R=1}^T\sum_{S=R+1}^{R+g(N,R)}\bp(A(R)A(S))+
\sum_{R=1}^T\sum_{S=R+g(N,R)+1}^T\bp(A(R)A(S))=:I+II.
\end{eqnarray*}

\noindent By (\ref{g2}) we have
\begin{eqnarray}
I&\leq&O(1)\sum_{R=1}^T\sum_{S=R+1}^{R+g(N,R)}
{\frac{2^R}{2^{S+g(N,S)}}}{\frac{1}{S+g( N,S)}}\leq\\
&\leq& O(1)\sum_{R=1}^T\frac{1}{(R+g(N,R))2^{g(N,R)}}
\sum_{j=1}^{g(N,R)}\frac{1}{2^{j}}
\leq\nonumber\\
&\leq& O(1)\sum_{R=1}^T \frac{1}{R2^{g(N,R)} }\leq
O(1)\sum_{R=1}^T \frac{1}{\lambda(N,R)}\leq O(1)(\log_N T)
\nonumber
\end{eqnarray}

\noindent and \begin{eqnarray} II&\leq&
O(1)\sum_{R=1}^T\sum_{S=R+g(N,R)+1}^T{\frac{1}{\lambda(N,R)}}
{\frac{1}{\lambda(N, S-R)}}\leq\\
&\leq& O(1)(\log_N T)^2.\nonumber
\end{eqnarray}

\noindent By (\ref{g1}) and (\ref{g2})

\begin{equation}\label{two}
\sum_{R=1}^T\sum_{S=R+1}^T\bp(A(R)A(S))\leq O(1)(\log_N T)^2.
\end{equation}

(\ref{g3}), (\ref{two}) and the Kochen-Stone Borel--Cantelli lemma
(see e.g. Spitzer \cite{SP}, page 317)  imply the first statement
with positive probability. Now to finish our proof  we need to
apply the zero-one law (again in a non-independent set up) as in
the proof of Theorem 3.2, observing   that  for any $n=1,2,\ldots$
$$\bp(A(R) \,\,{\rm i.o.}\mid X_1,X_2,
...X_n)=\bp(A(R) \,\,{\rm i.o.}).
$$
 $\Box$

\vspace{2ex}\noindent {\bf Proof of Theorem 3.5:} The proof goes
along the same line as the proof of Theorem 3.4.   The only point
which needs a little different approach is the the proof the
counterpart of Lemma 3.5. Namely, in the proof of this lemma we
need an upper bound for $1-p(R,R+1,S)$, which is equivalent of
getting a lower bound for D(R,S). Observe that in Lemma 2.4 we
have an asymptotic formula for $D(R,\infty).$  Now to get a lower
bound for $D(R,S)$ we need a less precise calculation (the
statement of Theorem  3.5 does not depend on $B,$ which was
important in Lemma 2.4). It is enough to observe that

$$U_i\geq C \exp\left(-\frac{B^*}{i^{\alpha} } \right)$$ with an
appropriate choice of $C$ and $B^*>B.$   After this observation,
with some tedious calculation somewhat similar to Lemma 2.4, we
get that
\begg D(R,S)\geq C
R^{\alpha}\left(1-\left(\frac{S}{R}\right)^{\alpha} \exp\left(
C_1(R^{1-\alpha}-S^{1-\alpha})\right) \right). \label{drs}\endd
 It is easy to see $$D(R,S)\geq C_2R^{\alpha}$$  if $S\geq
 R+R^{\alpha}/\log R.$ On the other hand, if $R<S < R+ R^{\alpha}/\log
 R$ then it can be seen that
 $$D(R,S)>C_3 (S-R)$$
and this is enough to carry through the argument in Lemma 3.5. We
omit the details. $\Box$


\bigskip\noindent
\begin{thebibliography}{99}
\bibitem{BRS} BR\'EZIS, H., ROSENKRANTZ, W. and SINGER, B.: An extension
of Khintchine's estimate for large deviations to a class of Markov
chains converging to a singular diffusion. {\em Comm. Pure Appl. Math.}
{\bf 24} (1971), 705--726.

\bibitem{CH} CHUNG, K.L.: {\em Markov Chains with Stationary
Transition Probabilities}, 2nd ed. Springer-Verlag, New York, 1967.

\bibitem{C-SD} COOLIN-SCHRIJNER, P. and VAN DOORN, E.A.:
Analysis of random walks using orthogonal polynomials.
{\em J. Comput. Appl. Math.} {\bf 99} (1998), 387--399.

\bibitem{CSFR06} CS\'AKI, E., F\"OLDES, A. and R\'EV\'ESZ, P.:
On the local time of the asymmetric Bernoulli walk (2006) (to be
submitted)

\bibitem{DE01} DETTE, H.: First return probabilities of birth and death
chains and associated orthogonal polynomials. {\em Proc. Amer. Math.
Soc.} {\bf 129} (2001), 1805--1815.
\bibitem{DO} DOOB, J.L.: {\em Stochastic Processes}. John Wiley
\& Sons, Inc. N.Y. (1953).
\bibitem{GA84} GALLARDO, L.: Comportement asymptotique des marches
al\'eatoires associ\'ees aux polyn\^omes de Gegenbauer. {\em Adv. Appl.
Probab.} {\bf 16} (1984), 293--323.

\bibitem{RG} GRADSTEYN, I.S. and RYZHIK I.M.: {\em Table of Integrals,
Series, and Products}. Academic Press, New York, 1980.

\bibitem{JLP} JAMES, N., LYONS, R. and PERES, Y.: A transient
Markov chain with finitely many cutpoints. {\em arXiv:0706.2013}

\bibitem{KMG} KARLIN, S. and McGREGOR, J.: Random walks.
{\em Illinois J. Math.} {\bf 3} (1959), 66--81.

\bibitem {KN81} KNIGHT, F.B.: {\em Essentials of Brownian Motion and
Diffusion}. Am. Math. Soc., Providence, R.I., 1981.

\bibitem{LA60} LAMPERTI, J.: Criterian for the recurrence or
transience of stochastic processes. I. {\em J. Math. Anal.
Appl.} {\bf 1} (1960), 314--330.

\bibitem{LA63} LAMPERTI, J.: A new class of probability limit theorems.
{\em J. Math. Mech.} {\bf 11} (1962), 749--772.

\bibitem{LA602} LAMPERTI, J.: Criteria for stochastic processes. II.
Passage-time moments. {\em J. Math. Anal. Appl.} {\bf 7} (1963),
127--145.

\bibitem{R05} R\'{E}V\'{E}SZ, P.: {\em Random Walk in Random and
Non-Random Environments}, 2nd ed. World Scientific, Singapore,
2005.

\bibitem{RO73} ROSENKRANTZ, W.A.: A method for computing the asymptotic
limit of a class of expected first passage times. {\em Ann. Probab.}
{\bf 1} (1973), 1035--1043.

\bibitem{SP} SPITZER, F.: {\em Principles of Random Walk}, Springer-Verlag, New York 1976.

\bibitem{SZ} SZ\'EKELY, G.J.: On the asymptotic properties of diffusion
processes. {\em Ann. Univ. Sci. Budapest, E\"otv\"os, Sect. Math.}
{\bf 17} (1974), 69--71.

\bibitem{VO90} VOIT, M.: A law of the iterated logarithm for a class of
polynomial hypergroups. {\em Monatshefte Math.} {\bf 109} (1990),
311--326.

\bibitem{VO92} VOIT, M.: Strong laws of large numbers for random walks
associated with a class of one-dimensional convolution structures. {\em
Monatshefte Math.} {\bf 113} (1992), 59--74.

\end{thebibliography}
\end{document}